\documentclass[12pt]{article}
\usepackage{amsfonts,amssymb}
\usepackage{url}
\newtheorem{defin}{Definition}[section]

\newtheorem{lem}[defin]{Lemma}
\newtheorem{lemma}[defin]{Lemma}

\newtheorem{thm}[defin]{Theorem}
\newtheorem{theorem}[defin]{Theorem}
\newtheorem{cor}[defin]{Corollary}

\newtheorem{claim}{Claim}
\newtheorem{con}[defin]{Conjecture}

\def\qed{\hbox{\kern1pt\vrule height6pt width4pt
depth1pt\kern1pt}\medskip}
\def\bproof{\par\noindent{\bf Proof.\enspace}\rm}

\newcommand{\real}{{\mathbb{R}}}
\newcommand{\complex}{{\mathbb{C}}}
\newcommand{\rat}{{\mathbb{Q}}}

\newcommand{\scrf}{{\mathcal{F}}}

\newcommand{\sm}{\setminus}

\newcommand{\eproof}{\hfill $\bullet$\\}

\newcommand{\Fiso}{{\cal F}_{iso}}

\title{Radically solvable graphs}
\date{1 July 2012}
\author{Bill Jackson \thanks{School of Mathematical Sciences, Queen Mary
University of London, Mile End Road, London E1 4NS, England. E-mail:
b.jackson@qmul.ac.uk} \and J.C. Owen \thanks{Siemens, Park House,
Cambridge CB3 0DU, England. E-mail: owen.john.ext@siemens.com}}
\begin{document}
\maketitle

\begin{abstract}
A 2-dimensional framework is a straight line realisation of a graph
in the Euclidean plane. It is radically solvable if the set of
vertex coordinates is contained in a radical extension of the field
of rationals extended by the squared edge lengths. We show that the
radical solvability of a generic framework depends only on its
underlying graph and characterise which planar graphs give rise to
radically solvable generic frameworks. We conjecture that our
characterisation extends to all graphs.
\end{abstract}


\section{Introduction}
Many systems of polynomial equations which are of practical interest
can be represented by a graph.
An important example occurs in computer aided design (CAD) when the
location of the geometric elements in a drawing such as points and
lines  (corresponding to vertices in the graph) are determined by
relationships between them such as tangency, coincidence and the
relative separations or angles between them (corresponding to edges
in the graph). The ability to solve such systems of equations
rapidly allows a design engineer to modify input parameters such as
the values for the separations or angles (collectively called
"dimensions" in a dimensioned drawing) and to realise a computer
model for many variants of a basic design \cite{Ow}. Most modern CAD
systems incorporate the ability to solve these so-called dimensional
constraint equations, see for example \cite{D3}.

A simple example of dimensional constraint equations is provided by
points in a plane with certain specified relative distances.
Both the equations and a particular solution can be represented by a framework
$(G,p)$ where $G$ is a graph and $p$ is a
vector comprising of all the
coordinates of the points. The graph $G$ has a vertex for each point and an
edge for each specified distance. Since the coordinates of the points are specified in
$(G,p)$ it is a simple matter to determine the relative distance corresponding
to any edge of $G$. The framework $(G,p)$ therefore represents both a system of polynomial equations and
a particular solution to these equations. We will call these equations the framework equations - they
correspond to the dimensional constraint equations referred to
above. In general the framework vector $p$ will be
just one of the many possible solutions to the framework equations.
(Estimates on the number of solutions have been obtained by several authors, see for example \cite{BS,JO,ST}.)


Efficient algorithms for solving the framework equations are
extremely useful. A particularly desirable case is when there are
only a finite number of solutions, and
these solutions can be
expressed as a sequence of square, or higher power, roots of
combinations of the squared edge distances. Such frameworks are said
to be quadratically solvable (or ruler-and-compass-constructible
\cite{Gao})
 and radically
solvable, respectively.
We will consider the problem of determining which generic frameworks
are quadratically or radically solvable.

The condition that the framework equations should have only finitely
many solutions is equivalent to the statement that the framework is
rigid. This property has been extensively studied and we refer the
reader to \cite{W} for an excellent survey of the area. Previous
work on quadratic/radical solvability \cite{Ow,OP} considered
generic frameworks which are minimally rigid i.e. cease to be rigid
when any edge is removed. A conjectured characterisation of
quadratically/radically solvable minimally rigid generic frameworks was given in
\cite{Ow} and this conjecture was verified for the special case when
the underlying graph is 3-connected and planar in \cite{OP}.

We will extend the study of quadratic and radical solvability to include generic
frameworks which are
rigid but not necessarily minimally rigid.
We first show in Lemma \ref{radgen} that  the quadratic or radical solvability of a generic
framework
depends only on the underlying graph.
This means that if a graph is
quadratically or radically solvable then there
will be a quadratic or radical solution to the corresponding system of framework equations for
any sufficiently general but consistent set of input distances.
We next consider globally rigid graphs
i.e. graphs for which every generic realisation is a unique
solution to the corresponding framework equations. We show in Theorem \ref{globallyrigid} that all such graphs are
quadratically solvable.

We develop a reduction scheme in Section \ref{3con} which shows how
the radical or quadratic solvability of a rigid graph is related to
the corresponding property for a derived graph which may be chosen
to be minimally rigid. We use this and the main result of \cite{OP}
to show in Theorem \ref{planar} that a rigid 3-connected planar
graph is radically solvable if and only if it is globally rigid.
This leads us to consider rigid graphs which are not 3-connected
i.e. graphs $G=(V,E)$ which can be separated into two subgraphs
$G=G_1 \cup G_2$  with $V(G_1) \cap V(G_2)=\{u,v\}$. We show in
Theorem \ref{2sep} that the radical or quadratic solvability of $G$
is determined by the corresponding property of $G_1+uv$ and $G_2+uv$
when $G_1$ and $G_2$ are both rigid, and of $G_1+uv$ and $G_2$ when
$G_1$ is not rigid and $G_2$ is minimally rigid. We use this
analysis to give a constructive definition for a family of
quadratically solvable graphs $\scrf$. We conjecture that every
radically solvable graph belongs to $\scrf$ and prove in Theorem
\ref{planar2con} that this holds for planar graphs.


\section{Definitions and Notation}\label{intro}
All graphs considered are finite and without loops or multiple
edges. Given a graph $G=(V,E)$ and two vertices $u,v\in V$ we use
$G+uv$ to denote the graph $(V, E\cup \{uv\})$.
  A {\em complex (real) realisation} of $G$ is
a map $p$ from $V$ to $\complex^2$ ($\real^2$). We also refer to the
ordered pair $(G,p)$ as a complex (real) {\em framework}. Although
we are mainly concerned with real frameworks, we will work with
complex frameworks since most of our methods require an
algebraically closed field and our results can still be applied to
the special case of real frameworks. Henceforth we assume that all
frameworks not specifically described as real, are complex. A
framework $(G,p)$ is {\em generic} if the set of all coordinates of
the points $p(v)$, $v\in V$, is algebraically independent over
$\rat$.

Let $V=\{v_1,v_2,\ldots,v_n\}$ and $E=\{e_1,e_2,\ldots,e_m\}$. Given
a realisation $(G,p)$ of $G$ in $\complex^2$ and two vertices
$v_i,v_j\in V$ with $p(v_i)-p(v_j)=(a,b)$ put $d_p(v_i,v_j)=a^2+b^2$
and $d_p(e)=d_p(v_i,v_j)$ when $e=v_iv_j\in E$.  Two realistions
$(G,p)$ and $(G,q)$ are {\em equivalent} if $d_p(e)=d_q(e)$ for all
$e\in E$, and are {\em congruent} if $d_p(v_i,v_j)=d_q(v_i,v_j)$ for
all $v_i,v_j\in V$. The {\em rigidity map} $d_G:\complex^{2n}\to
\complex^m$ is defined by putting
$d_G(p)=(d_p(e_1),d_p(e_2),\ldots,d_p(e_m))$.
Thus $(G,p)$ and
$(G,q)$ are {equivalent} if and only if $d_G(p)=d_G(q)$. Note that,
if $(G,p)$ and $(G,q)$ are real frameworks, then they are equivalent
if and only if they have the same edge lengths and they are
congruent if and only if we can transform one to the other by
applying an isometry of $\real^2$ i.e. a translation, rotation or
reflection of the Euclidean plane.

A framework  is {\em globally rigid} if  all equivalent frameworks
are congruent to it. A real framework $(G,p)$ is {\em rigid} if
there exists an $\epsilon>0$ such that every real framework $(G,q)$
which is equivalent to $(G,p)$ and satisfies
$d(p(v)-q(v))=\|p(v)-q(v)\|^2<\epsilon$ for all $v\in V$, is
congruent to $(G,p)$.\footnote{Equivalently, a real framework
$(G,p)$ is rigid if every continuous motion of the points $p(v)$,
$v\in V$, in $\real^2$ which preserves the edge distances results in
a framework which is congruent to $(G,p)$.} It is known that both
the rigidity and the global rigidity of a generic framework depend
only on its underlying graph. We say that a graph $G$ is {\em rigid}
if some, or equivalently every, generic real realisation of $G$ is
rigid, and that $G$ is {\em globally rigid} if some, or equivalently
every, generic realisation of $G$ is globally rigid.

Let $K,L$ be fields with $K\subseteq L$. Then $L$ is a {\em radical
extension} of $K$ if there exist fields $K=K_1\subset
K_2\subset\ldots\subset K_t=L$ such that for all $1\leq i<t$,
$K_{i+1}=K_i(x_i)$ with $x_i^{n_i}\in K_i$ for some natural number
$n_i$. The field $L$ is a {\em quadratic extension} of $K$ if it is a
radical extension with $n_i=2$ for all $1\leq i< t$. We say that
$L:K$ is {\em radically solvable}, respectively {\em quadratically
solvable},  if $L$ is contained in a radical, respectively
quadratic, extension of $K$. A realisation $(G,p)$ of a rigid graph
$G$ is {\em radically solvable}, respectively {\em quadratically
solvable}, if there exists a congruent realisation $(G,q)$ such that
$\rat(q):\rat(d_G(q))$ is
radically, respectively quadratically, solvable.

\section{Field extensions and algebraic varieties}

The above definitions of radically and quadratically solvable field
extensions immediately imply the following result.

\begin{lem}\label{basic} Let $K\subseteq L\subseteq M$ be  fields. Then
$M:K$ is radically, respectively quadratically, solvable if and only
if $M:L$ and $L:K$ are both radically, respectively quadratically,
solvable.
\end{lem}

We next recall some definitions and results from Galois theory.
We adopt the notation of \cite{Stew} and refer the reader to this
text for further information on the subject.

Given a field extension $L:K$ we use $[L:K]$ to denote the {\em
degree} of the extension i.e. the dimension of $L$ as a vector space
over $K$. The extension is {\em finite} if it has finite degree. It
is {\em normal} if $L$ is the splitting field of some polynomial
over $K$. When $L:K$ is finite, a {\em normal closure of $L$ over
$K$} is a field $N$ such that $L\subseteq N$, $N:K$ is normal, and,
subject to these conditions, $N$ is minimal with respect to
inclusion. It is known that normal closures exist, are finite, and
are unique up to isomorphism, see \cite[Theorem 11.6]{Stew}. The
{\em Galois group}  $\Gamma(L:K)$ is the group of all automorphisms
of $L$ which leave $K$ fixed. Galois theory gives us the following
close relationship between radically/quadratically solvable
extensions and Galois groups, see  \cite[Theorems 15.3,
18.18]{Stew}.\footnote{The references to \cite{Stew} in this section
only give results on radically solvable extensions, but similar
proofs work for the special case of quadratically solvable
extensions.}

\begin{thm}\label{Galois} Let $K$ be a field of characteristic zero and
$N:K$ be a normal field extension. Then\\
(a) $N:K$ is radically solvable  if and only if $\Gamma(N:K)$ is a solvable group.\\
(b) $N:K$ is quadratically solvable if and only if
$\Gamma(N:K)$ is a $2$-group.
\end{thm}

Our next result allows us to decide whether a field extension $L:K$
is
 radically, respectively quadratically, solvable by
applying Theorem \ref{Galois} to its normal closure. This will be
used to show that the radical or quadratic solvability of a generic
framework depends only on its underlying graph.

\begin{figure}
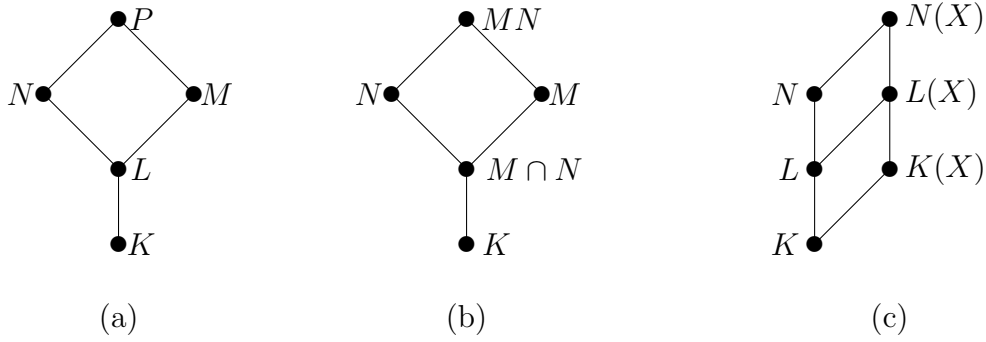

\input{radfig1a.tec}
\input{radfig1d.tec}
\input{radfig1h.tec}
\caption{The field extensions of Lemmas \ref{normal}, \ref{prod} and
\ref{indext} } \label{fig1}
\end{figure}

\begin{lem}\label{normal} Let $K$ be a field of characteristic zero
and $L:K$ be a finite field extension. Let $N$ be a normal closure
of $L$ over $K$. Then $L:K$ is radically, respectively
quadratically, solvable if and only if $N:K$ is radically,
respectively quadratically, solvable.
\end{lem}
\bproof Sufficiency follows from Lemma \ref{basic}. To prove
necessity we assume that $L:K$ is radically, respectively
quadratically, solvable. Then $L$ is contained in a radical,
respectively quadratic, extension $M$ of $K$. Let $P$ be a normal
closure of $M$ over $K$, see Figure \ref{fig1}(a). Since $L\subseteq
M$ and normal closures are unique up to isomorphism, we may suppose
that $N\subseteq P$. Since $M$ is a radical, respectively quadratic,
extension of $K$ and $P$ is a normal closure of $M$ over $K$,
\cite[Lemma 15.4]{Stew} implies that $P$ is also a radical,
respectively quadratic, extension of $K$. Since $N\subseteq P$,
$N:K$ is radically, respectively quadratically, solvable. \eproof

Suppose $M,N$ are field extensions of a field $K$ which are both
contained in a common extension $P$ of $K$. Then $MN$ denotes the
smallest subfield of $P$ which contains both $M$ and $N$. We will
need the following result from Galois Theory, see for example
\cite[Proposition 3.18]{M}.

\begin{lem}\label{prod} Let $K$ be a field of characteristic zero
and $M,N$ be field extensions of $K$ which are both contained in a
common extension of $K$. Suppose that $N$ is a normal extension of
$K$. Then $MN:M$ and $N:M\cap N$ are normal extensions, and
$\Gamma(MN:M)$ and $\Gamma(N:M\cap N)$ are isomorphic groups.
\end{lem}

Given a field $K$ we use $K[X_1,X_2,\ldots,X_n]$ to denote the ring
of polynomials in the indeterminates $X_1,X_2,\ldots,X_n$ with
coefficients in $K$ and $K(X_1,X_2,\ldots,X_n)$ to denote its field
of fractions.

\begin{lem}\label{indext}
Let $L:K$ be a finite field extension with $\rat\subseteq K\subseteq
L\subset\complex$, and $N$ be the normal closure of $L:K$ in
$\complex$. Let $X=(X_1,X_2,\ldots,X_n)$ be a vector of
indeterminates. Then $N(X)$ is a normal closure of $L(X)$ over
$K(X)$ and $\Gamma(N:K)$ is isomorphic to $\Gamma(N(X):K(X))$.
Furthermore, $L:K$ is radically, respectively quadratically,
solvable  if and only if $L(X):K(X)$ is radically, respectively
quadratically, solvable.
\end{lem}
\bproof
Let $a_1,a_2,\ldots,a_m$ be a basis for $L:K$, $f_i$ be the minimum
polynomial of $a_i$ over $K$, $R_i$ be the set of all complex roots
of $f_i$, and $R=\bigcup_{i=1}^mR_i$. Then $N=L(R)$. Since
$X_1,X_2,\ldots,X_n$ are indeterminates, $a_1,a_2,\ldots,a_m$ is
also a basis for $L(X):K(X)$ and $f_i$ is the minimum polynomial of
$a_i$ over $K(X)$. Thus $L(R)(X)=N(X)$ is a normal closure of
$L(X):K(X)$. We now apply Lemma \ref{prod} with $M=K(X)$. We have
$NK(X)=N(X)$ and $N\cap K(X)=K$. Hence $\Gamma(N:K)$ is isomorphic
to $\Gamma(N(X):K(X))$.

The final part of the lemma now follows from Theorem \ref{Galois}
and Lemmas \ref{normal} and \ref{prod}. \eproof

 Our next result is an application of the previous lemmas.
We will use it to determine whether generic realisations of graphs
with small separating sets of vertices are radically or
quadratically solvable.

\begin{lem}\label{sep} Suppose that $X=(X_1,X_2,\ldots,X_r)$, $Y=(Y_1,Y_2,\ldots,Y_s)$ and
$Z=(Z_1,Z_2,\ldots,Z_t)$ are vectors of indeterminates,
$f=(f_1,f_2,\ldots,f_m)\in \rat[X,Y]^m$ and
$g=(g_1,g_2,\ldots,g_n)\in \rat[Y,Z]^n$, and $\rat(X,Y,Z)$ is
a  finite extension of $\rat(f,g)$. Then $\rat(X,Y,Z):\rat(f,g)$ is
  radically, respectively quadratically,  solvable  if and only if $\rat(f,Y,Z):\rat(f,g)$
  and $\rat(X,Y):\rat(f,Y)$ are  both radically, respectively quadratically,  solvable.
\end{lem}
\bproof This follows from Lemma \ref{basic} (which tells us that
$\rat(X,Y,Z):\rat(f,g)$ is  radically, respectively quadratically,
solvable  if and only if $\rat(f,Y,Z):\rat(f,g)$ and
$\rat(X,Y,Z):\rat(f,Y,Z)$ are both radically, respectively
quadratically,  solvable) and Lemma \ref{indext} (which tells us
that $\rat(X,Y,Z):\rat(f,Y,Z)$ is radically, respectively
quadratically,  solvable if and only if  $\rat(X,Y):\rat(f,Y)$ is
radically, respectively quadratically,  solvable). \eproof

Our final result of this section concerns algebraic varieties. We
will use it to show, amongst other things, that globally rigid graphs are quadratically
solvable.

\begin{lem} \label{extend}
Let $K$ be a field with $\rat \subseteq K \subset \complex$, and let
$S=\{f_1,\dots,f_m\} \subset K[X]$, where $X=(X_1,\dots,X_n)$ is a
vector of indeterminates. Let $I \subset K[X]$ be the ideal
generated by $S$ and $W=\{x\in \complex^n\,:\,f_i(x)=0 \mbox{ for
all $1\leq i\leq m$}\}$. Let $I_1=I\cap K[X_1]$.
Then $I_1$ is an ideal of $K[X_1]$ and is generated by a single
polynomial $h_1\in K[X_1]$. Furthermore, if $W$ is non-empty and
finite and $h_1(a)=0$ for some $a\in \complex$, then there exists an
$x=(x_1,x_2,\ldots,x_n)\in W$ such that $x_1=a$.
\end{lem}
\bproof It is easy to see that $I_1$ is an ideal of $K[X_1]$. It is
generated by a single polynomial since  $K[X_1]$ is a principal
ideal domain.  The final part of the lemma follows from the work of
Kalkbrener \cite{K}. We include an outline of his proof for
completeness. Let $I_s=I\cap K[X_1,\ldots,X_s]$ for all $1\leq s\leq
n$ and $W_s=\{x\in \complex^s\,:\,f(x)=0 \mbox{ for all $f\in
I_s$}\}$. It will suffice to show that for all $2\leq s\leq n$ and
all $b\in W_{s-1}$, there exists a $c\in \complex$ such that
$(b,c)\in W_{s}$. Let $I_s(b)=\{f(b,X_s)\,:\,f\in I_s\}$. Then
$I_s(b)$ is an ideal of $K(b)[X_s]$. It follows from \cite[Theorems
3]{K} that
$I_s(b)$ is generated by a non-constant polynomial $h_s$ in
$K(b)[X_s]$. We can now choose a $c\in \complex$ with $h_s(c)=0$ to
obtain the required element $(b,c)\in W_s$.
\eproof

\section{Standard positions}
Given a generic framework $(G,p)$ it will be useful to identify a
particular congruent framework $(G,q)$ with the property that
$(G,p)$ is  radically, or quadratically, solvable if and only if
$\rat(q)$ is contained in a radical, or quadratic, extension of
$\rat(d_G(q))$. The following result will enable us to do this.

\begin{lem}\label{standard}
Suppose that $(G,p)$ is a generic realisation of a graph $G=(V,E)$ where
$V=\{v_1,v_2,\ldots,v_n\}$ and
$n\geq 3$.
Then
there are exactly four realisations $(G,q_j)$, $1\leq j\leq 4$, which are congruent to $(G,p)$
and have $q_j(v_1)=(0,0)$ and $q_j(v_2)=(0,z)$ for some $z\in \complex$.
Furthermore, we have $\rat(q_i)=\rat(q_j)$ for all $1\leq i<j\leq 4$. \end{lem}
\bproof
The assertion that there are exactly four such realisations $(G,q_i)$
is a special case of \cite[Corollary 5.3]{JO}. The assertion that $\rat(q_i)=\rat(q_j)$ follows
from the fact that we can order the $q_j$ such that, if $q_1(v_i)=(x_i,y_i)$ for all $v_i\in V$,
then $q_2(v_i)=(-x_i,y_i)$, $q_3(v_i)=(x_i,-y_i)$ and $q_4(v_i)=(-x_i,-y_i)$ for all $v_i\in V$.
\eproof

Given a graph $G$ and vertices $v_1,v_2$ of $G$, we say that a
realisation $(G,q)$ of $G$ is in {\em standard position with respect
$(v_1,v_2)$} if $q(v_1)=(0,0)$ and $q(v_2)=(0,z)$ for some $z\in
\complex$, and is {\em quasi-generic} if it is congruent to a
generic realisation of $G$.

\begin{lem}\label{algind}
Suppose that $(G,p)$ is a quasi-generic realisation of a rigid graph
$G=(V,E)$ where $V=\{v_1,v_2,\ldots,v_n\}$ and $p(v_i)=(x_i,y_i)$
for $1\leq i\leq n$. Suppose further that $(G,p)$ is in standard
position with respect to $(v_1,v_2)$, i.e. $x_1=y_1=x_2=0$. Then
$\{y_2,x_3,y_3,\ldots,y_n\}$ is algebraically independent over
$\rat$.
\end{lem}
\bproof This follows immediately from \cite[Lemma 5.4]{JO}.
\eproof

\begin{lem}\label{radical}
Suppose that $(G,p)$ is a generic realisation of a rigid graph
$G=(V,E)$ where $V=\{v_1,v_2,\ldots,v_n\}$ and $n\geq 3$, and
$(G,q)$ is a congruent realisation in standard position with respect
to $(v_1,v_2)$. Then $(G,p)$ is radically, respectively
quadratically, solvable if and only if $\rat(q):\rat(d_G(q))$ is
radically, respectively quadratically, solvable.
\end{lem}
\bproof Sufficiency follows immediately from the definition of
radically, respectively quadratically, solvable frameworks. To prove necessity
we suppose that $(G,p)$ is radically, respectively quadratically,
solvable. Replacing $(G,p)$ by a congruent framework if necessary, we
may assume that $\rat(p)$ is itself contained in a radical,
respectively quadratic,  extension $L$ of $\rat(d_G(p))$. We
can construct a framework $(G,q)$ satisfying the hypotheses of the
lemma by putting $\tilde q(v_i)= p(v_i)-p(v_1)$ for all $v_i\in V$,
and
$$q(v_i)=\left(
\begin{array}{rr}
y/d_0 & -x/d_0\\
x/d_0 & y/d_0
\end{array}
\right)\tilde q(v_i)$$ for all $v_i\in V(G)$, where $\tilde
q(v_2)=(x,y)$ and $d_0^2 = x^2+y^2$. By Lemma \ref{standard}, it
will suffice to show that for this $q$, $\rat(q)$ is contained in  a
radical, respectively quadratic, extension of
$\rat(d_G(q))$. Let $K=\rat(p,d_0)$. The definitions of $\tilde q$
and $q$ imply that $\rat(\tilde q)\subseteq \rat(p)$ and hence that
$\rat(q)\subseteq K$. We have $[K:\rat(p)]\leq 2$ since $d_0^2 =
x^2+y^2$ and $x,y\in \rat(p)$. Hence $L(d_0)$ is a radical,
respectively quadratic,  extension of $\rat(d_G(p))$ which contains $K$.
Since $\rat(q)\subseteq K$ and $d_G(p)=d_G(q)$,
$\rat(q):\rat(d_G(q))$ is radically, respectively quadratically,
solvable. \eproof

\section{Quadratically and radically solvable graphs}

We first show that a quasi-generic realisation of a rigid graph
gives rise to a finite field extension when it is in standard
position.

\begin{lem}\label{finite}
Suppose that $G=(V,E)$ is a rigid graph and that $(G,p)$ is a quasi-generic realisation of $G$ in standard position with respect to
two vertices $v_1,v_2\in V$. Then $\rat(p):\rat(d_G(p))$ is a finite field extension.
\end{lem}
\bproof
It is easy to see that $\rat(d_G(p))\subseteq \rat(p)$. By \cite[Lemma 5.4]{JO},  $\rat(p)$ and $\rat(d_G(p))$
have the same algebraic closure. This implies that each coordinate of $p$ is a root of a polynomial
with coefficients in $\rat(d_G(p))$ and hence  $[\rat(p):\rat(d_G(p))]$ is finite.
\eproof

  We next show that radical and quadratic solvability are generic properties of frameworks i.e.
  they depend only on the underlying graph when the given realisation is
  generic.

\begin{lem}\label{radgen}
Suppose $(G,p)$ and $(G,p')$ are generic realisations of a rigid
graph $G=(V,E)$. Then $(G,p)$ is radically, respectively
quadratically, solvable if and only if  $(G,p')$ is radically,
respectively quadratically, solvable.
\end{lem}
\bproof Let $V=\{v_1,v_2,\ldots,v_n\}$ and
$E=\{e_1,e_2,\ldots,e_m\}$. Let $(G,q)$ and $(G,q^\prime)$ be two
frameworks in standard position with respect to $(v_1,v_2)$ which
are congruent to $(G,p)$ and $(G,p^\prime)$, respectively. Put
$q(v_i)=(x_{2i-1},x_{2i})$ and $q'(v_i)=(x'_{2i-1},x'_{2i})$ for
$1\leq i\leq n$. We associate a pair of indeterminates
$(X_{2i-1},X_{2i})$ with each vertex $v_i\in V$, putting
$X_1=X_2=X_3=0$ to represent a framework in standard position. Let
$X=(X_4,X_5,\ldots,X_{2n})$ and $D_G(X)=(f_1,f_2,\ldots,f_m)$ where
$f_i=(X_{2j-1}-X_{2k-1})^2+(X_{2j}-X_{2k})^2$ when $e_i=v_jv_k$.
Since $(G,q)$ and $(G,q')$ are quasi-generic, Lemma \ref{algind}
implies that $\{x_3,x_4,\ldots,x_{2n}\}$ and
$\{x'_3,x'_4,\ldots,x'_{2n}\}$ are both algebraically independent
over $\rat$. Hence $\rat(q):\rat(d_G(q))$ and
$\rat(q'):\rat(d_G(q'))$ are both isomorphic to
$\rat(X):\rat(D_G(X))$.\footnote{Two field extensions $L:K$ and
$L':K'$ are isomorphic if there exists a field isomorphism from $L$
to $L'$ which maps $K$ onto $K'$.}

Let $N_q$, $N_{q'}$ and $N_X$ be normal closures of
$\rat(q):\rat(d_G(q))$, $\rat(q'):\rat(d_G(q'))$ and
$\rat(X):\rat(D_G(X))$, respectively. Then $N_q:\rat(d_G(q))$ and
$N_{q'}:\rat(d_G(q'))$ are both isomorphic to $N_X:\rat(D_G(X))$ and
hence are isomorphic to each other. It follows that
$\Gamma(N_q:\rat(d_G(q)))$ and $\Gamma(N_{q'}:\rat(d_G(q')))$ are
isomorphic groups. The lemma now follows by applying Theorem
\ref{Galois} and Lemma \ref{normal}.
\eproof

This result allows us to  define a rigid graph to be {\em
radically}, respectively {\em quadratically}, {\em solvable} if some
(or equivalently every) generic realisation of $G$ is {radically},
respectively quadratically, {solvable}. Lemmas \ref{radical} and
\ref{radgen} imply that this definition agrees with the one given
for the radical and quadratic solvability of minimally rigid graphs
in \cite[Definition 3.1]{OP}.

\section{Globally rigid graphs}

Two vertices $v_i,v_j$ of a rigid graph $G$ are {\em globally
linked} if for each generic realisation $(G,p)$ and every equivalent
realisation $(G,q)$ we have $d_p(v_i,v_j)=d_q(v_i,v_j)$.

\begin{lem} \label{globally_b}
Let $(G,p)$  be a quasi-generic realisation of a rigid graph
$G=(V,E)$ with $V=\{v_1,v_2,\ldots,v_n\}$ and
$E=\{e_1,e_2,\ldots,e_m\}$. Suppose that $v_a,v_b \in V$ are
globally linked in $G$. Then $d_p(v_a,v_b) \in \rat(d_G(p))$.
\end{lem}
\bproof  We may suppose that $(G,p)$ is in standard position with
respect to $v_1,v_2$. Let $K=\rat(d_G(p))$. We again associate a
pair of indeterminates $(X_{2i-1},X_{2i})$ with each vertex $v_i\in
V$, putting $X_1=X_2=X_3=0$ to represent a framework in standard
position. Let
$f_i=(X_{2j-1}-X_{2k-1})^2+(X_{2j}-X_{2k})^2-d((p(v_j)-p(v_k))$ for
each $e_i=v_jv_k\in E$. We introduce a new indeterminate $X_{2n+1}$
which represents the `distance' between $v_a$ and $v_b$ and put
$f_{m+1}=X_{2n+1}-(X_{2a-1}-X_{2b-1})^2-(X_{2a}-X_{2b})^2$. Let
$X=(X_4,X_5,\ldots,X_{2n+1})$. Let  $I$ be the ideal of $K[X]$
generated by the polynomials $f_1,f_2,\ldots,f_{m+1}$ and let
$I_{2n+1}=I\cap K[X_{2n+1}]$. Then $I_{2n+1}$ is generated by a
single polynomial $h_{2n+1}\in K[X_{2n+1}]$, and every zero of
$h_{2n+1}$ in $\complex$ extends to a zero of $I$ in
$\complex^{2n+1}$ by Lemma \ref{extend}. Since $v_a,v_b$ are
globally linked in $G$, $d_p(v_a,v_b)$ must be the unique zero of
$h_{2n+1}$. Thus $h_{2n+1}=(X_{2n+1}-d_p(v_a,v_b))^t$ for some
positive integer $t$. Since $h_{2n+1}\in K[X_{2n+1}]$ this implies
that $d_p(v_a,v_b)\in K$. \eproof

\begin{thm} \label{globallyrigid}
Every globally rigid graph is quadratically solvable.
\end{thm}
\bproof Let $G=(V,E)$ be a  globally rigid graph with
$V=\{v_1,v_2,\ldots,v_n\}$ and $E=\{e_1,e_2,\ldots,e_m\}$. Let
$(G,p)$ be quasi-generic realisation of $G$ which is in standard
position with $p(v_1)=(0,0)$ and $p(v_2)=(0,y_2)$. Let
$K=\rat(d_G(p))$ and $K_1=K(y_2)$. Since $y_2$ satisfies the
quadratic  equation $y_2^2-d_p(v_1,v_2)=0$, and since
$d_p(v_1,v_2) \in K$ by Lemma \ref{globally_b},
we have $[K_1:K] \leq 2$.
Let $p(v_i)=(x_i,y_i)$ for all $3\leq i\leq n$.
Then $x_i^2+y_i^2=d_p(v_i,v_0)$ and
$x_i^2+(y_i-y_2)^2=d_p(v_i,v_1)$.
  Since $G$ is globally rigid,
$v_i$ is globally linked to both $v_1$ and $v_2$ in $G$ and hence,
by Lemma \ref{globally_b}, $\{d_p(v_i,v_0), d_p(v_i,v_1)\} \subset
K$. This implies that $y_i\in K_1$ and $x_i^2\in K_1$. Since this
holds for all $3\leq i\leq n$, $(G,p)$ is quadratically solvable.
\eproof

\section{3-connected graphs}\label{3con}

A graph $G=(V,E)$ is {\em $k$-connected} if $|V|\geq k+1$ and $G-U$
is connected for all $U\subseteq V$ with $|U|<k$.
We conjecture that a 3-connected graph is radically (or
quadratically) solvable if and only if it is globally rigid. We will
verify this conjecture for planar graphs. In addition we show that
our conjecture is equivalent to an old conjecture of the second
author (that no 3-connected minimally rigid graph is radically
solvable). We will use the following lemma which tells us that the
radical, respectively quadratic, solvability of a rigid graph is
preserved by the operation of replacing a subgraph by a  radically,
respectively quadratically, solvable rigid subgraph. (In our
application the new subgraph will be minimally rigid.)

\begin{lem} \label{replace}
Let $H_0,H_1,H_2$ be graphs with $V(H_0)\cap V(H_1)=V(H_0)\cap
V(H_2)=V(H_1)\cap V(H_2)=U$, $|U|\geq 2$, and $E(H_0)\cap
E(H_1)=E(H_0)\cap E(H_2)=\emptyset$. Let $G_1=H_0\cup H_1$ and
$G_2=H_0\cup H_2$.
Suppose that $G_1$ and $H_2$ are both rigid. Then\\
(a) $G_2$ is rigid.\\
(b) If $G_1$ and $H_2$ are both radically, respectively
quadratically, solvable then $G_2$ is radically, respectively
quadratically, solvable.
\end{lem}
\bproof
Choose $v_1,v_2\in U$ and let $(G_1\cup G_2,p)$ be a quasi-generic
{\em real} realisation of $G_1\cup G_2$ with $p(v_1)=(0,0)$ and
$p(v_2)=(0,y)$ for some $y\in \real$. Let $V_i=V(H_i)\sm U$ for
$0\leq i\leq 2$.

Suppose that $G_2$ is not rigid. Since $H_2$ is rigid, there exists
a non-zero infinitesimal motion $z_2$ of $(G_2,p)$ in $\real^2$
which keeps $H_2$ fixed. Then $z_1:V(G_1)\to \real^2$ by
$z_1(v)=(0,0)$ for $v\in V(H_1)$ and $z_1(v)=z_2(v)$ for $v\in
V(H_0)$  is a non-zero infinitesimal motion of $G_1$ which keeps
$H_1$ fixed. This contradicts the hypothesis that $G_1$ is rigid and
completes the proof of (a).

Suppose that $G_1$ and $H_2$ are both radically, respectively
quadratically, solvable. The first assumption implies that
$\rat(p|_{V_0},p|_{U},p|_{V_1})$ is
a radically,
respectively quadratically, solvable extension of
$\rat(d_{H_0}(p),d_{H_1}(p))$. Since the components of $(p|_{V_0},y,p|_{U\sm
\{v_1,v_2\}},p|_{V_1})$ are algebraically independent over $\rat$ we
may treat them
as if they were indeterminates and apply
Lemma \ref{sep} with $X=p|_{V_0}$, $Y=(y,p|_{U\sm\{v_1,v_2\}})$,
$Z=p|_{V_1}$, $f=d_{H_0}(p)$, and $g=d_{H_1}(p)$ to deduce that
$\rat(p|_{V_0},p|_U)$ is
a radically, respectively
quadratically, solvable extension of $\rat(d_{H_0}(p),p|_U)$. We
also have $\rat(p|_U,p|_{V_2})$
is a radically, respectively quadratically, solvable extension of
$\rat(d_{H_2}(p))$ by the second assumption. Hence
$\rat(d_{H_0}(p),p|_U,p|_{V_2})$ is
a radically, respectively quadratically, solvable extension of
$\rat(d_{H_0}(p),d_{H_2}(p))$. Since the components of $(p|_{V_0},y,p|_{U\sm
\{v_1,v_2\}},p|_{V_2})$ are algebraically independent over $\rat$, we
may apply Lemma \ref{sep}, with $X=p|_{V_0}$,
$Y=(y,p|_{U\sm\{v_1,v_2\}})$, $Z=p|_{V_2}$, $f=d_{H_0}(p)$, and
$g=d_{H_2}(p)$, to deduce  that $\rat(p|_{V_0},p|_{U},p|_{V_2})$ is
a radically, respectively quadratically, solvable extension of
$\rat(d_{H_0}(p),d_{H_2}(p))$. Thus $G_2$ is radically, respectively
quadratically, solvable and (b) holds. \eproof

We also need a result on graph connectivity due to W. Mader.

\begin{lem}\cite[Satz 1]{Mad} \label{mad}
Let $G$ be a $k$-connected graph and $C$ be a cycle in $G$ such that
each vertex of $C$ has degree at least $k+1$ in $G$. Then $G-e$ is
$k$-connected for some $e\in E(C)$.
\end{lem}

For $n\geq 4$, the {\em wheel on $n$ vertices} is the graph
$W=(V,E)$ with $V=\{v,u_1,\ldots,u_{n-1}\}$ and
$E=\{vu_1,vu_2,\ldots,vu_{n-1}\}\cup\{u_1u_2,u_2u_3,\ldots,u_{n-1}u_1\}$.
We refer to the cycle $C=u_1u_2\ldots u_{n-1}u_1$ as the {\em rim}
of $W$, and to the vertices of $C$ as the {\em rim vertices} of $W$.

\begin{lem} \label{3conreplace}
Let $H_0,H_1$ be graphs with $V(H_0)\cap V(H_1)=U$, $|U|\geq 3$, and
$E(H_0)\cap E(H_1)=\emptyset$. Let $H_2$ be a wheel with $U$ as its
set of rim vertices, $V(H_0)\cap V(H_2)=U$ and $E(H_0)\cap
E(H_2)=\emptyset$. Put $G_1=H_0\cup H_1$ and $G_2=H_0\cup H_2$.
Suppose that $G_1$ is $3$-connected and that each vertex of $U$ has
degree at least four in $G_2$. Then $G_2-e$ is $3$-connected for
some edge $e$ on the rim of $H_2$. Furthermore, if $G_1$ is planar
and $H_1$ is connected, then we may choose $H_2$ in such a way such
that $G_2-e$ is planar and $3$-connected.
\end{lem}
\bproof We first show that $G_2$ is 3-connected. Suppose not. Then
$G_2-T$ is disconnected for some $T\subseteq V(G_2)$ with $|T|\leq
2$. Since $H_2$ is 3-connected,
$H_2-T$ is connected. Hence $H_2-T$ is contained in a single
connected component of $G_2-T$. This implies that $G_1-(T\cap V(G_1))$ is
disconnected and contradicts the hypothesis that $G_1$ is
3-connected.

We may now use Lemma \ref{mad} and the hypothesis that each vertex
of $U$ has degree at least four in $G_2$ to deduce that $G_2-e$ is
3-connected for some edge $e$ of $C$.

Finally, we suppose that $G_1$ is planar and $H_1$ is connected.
Then the vertices of $U$ must lie on the same face $F$ of
$G-(V(H_1)-U)$. If we choose $H_2$ such that, in the above
definition of a wheel, the rim vertices $u_1,u_2,\ldots,u_{n-1}$
occur in this order around $F$, then the resulting $G_2$ will be
planar. \eproof


\begin{lem}\label{wheel} Let $G$ be obtained by deleting an edge from the rim of
a wheel on $n\geq 4$ vertices. Then $G$ is both minimally rigid and
quadratically solvable.
\end{lem}
\bproof
It is easy to check that $G$ can be obtained from $K_3$ by
recursively adding vertices of degree two. The lemma now follows
since $K_3$ is minimally rigid and quadratically solvable, and the
operation of adding a vertex of degree two is known to preserve the
properties of being minimally rigid, see \cite{W}, and quadratically
solvable \cite{Ow}. \eproof

A graph $G=(V,E)$ is {\em redundantly rigid} if $G-e$ is rigid for
all $e\in E$. A {\em non-trivial redundantly rigid component} of $G$
is a maximal redundantly rigid subgraph of $G$. Edges $e$ of $G$
such that $G-e$ is not rigid belong to no redundantly rigid
subgraphs of $G$. We consider the subgraph consisting of such an
edge $e$ and its end-vertices to be a {\em trivial redundantly rigid
component}. Thus  $G$ is minimally rigid if and only if all its
redundantly rigid components are trivial and, when $|V|\geq 3$, $G$
is redundantly rigid if and only if it has exactly one redundantly
rigid component.

We can now characterise quadratic solvability in 3-connected planar
graphs. We use the fact that a rigid graph $G=(V,E)$ is minimally
rigid if and only if $|E|=2|V|-3$, see \cite{W}.

\begin{thm}\label{planar}
Let $G=(V,E)$ be a rigid $3$-connected planar graph. Then the
following
statements are equivalent.\\
(a) $G$ is quadratically solvable.\\
(b) $G$ is radically solvable.\\
(c) $G$ is redundantly rigid.\\
(d) $G$ is globally rigid.
\end{thm}
\bproof If $G$ is redundantly rigid then $G$ is globally rigid by
\cite{JJ} and hence is quadratically solvable by Theorem
\ref{globallyrigid}. Hence (c) implies (d) and (d) implies (a).
Clearly (a) implies (b). It remains to show that (b) implies (c). We
will prove the contrapositive.

Suppose that $G$ is not redundantly rigid. We show by induction on
$|E|-2|V|+3$ that $G$ is not quadratically solvable. Since $G$ is
rigid we have $|E|- 2|V|+3\geq 0$. If equality holds then $G$ is
minimally rigid and \cite{OP} implies that $G$ is not radically
solvable. Hence we may suppose that $|E|> 2|V|-3$. Then some
redundantly rigid component $H_1=(V_1,E_1)$ of $G$ is non-trivial.
Let $U$ be the set of vertices of $H_1$ which are incident to edges
of $E\sm E_1$ and put $H_0=(G-E_1)-(V_1\sm U)$. By Lemma
\ref{3conreplace}, we can choose a wheel $W$ with rim vertices $U$
and an edge $e$ on the rim of $W$ such that $G'=H_0\cup (W-e)$ is
3-connected and planar. Lemmas \ref{replace}(a) and \ref{wheel}
imply that $G'$ is rigid. Since $G'$ is not redundantly rigid and
$|V(G')|-2|E(G')|+3<|E|- 2|V|+3$, we may apply induction to deduce
that $G'$ is not radically solvable. Lemmas \ref{replace}(b) and
\ref{wheel} now imply that $G$ is not radically  solvable. \eproof

We conjecture that the planarity condition can be removed from
Theorem \ref{planar}.

\begin{con}\label{gen}
Let $G=(V,E)$ be a rigid $3$-connected graph. Then the following
statements are equivalent.\\
(a) $G$ is quadratically solvable.\\
(b) $G$ is radically solvable.\\
(c) $G$ is redundantly rigid.\\
(d) $G$ is globally rigid.
\end{con}

We may use the proof technique of Theorem \ref{planar} to reduce
this conjecture to the special case when $G$ is minimally rigid.
This special case was suggested over twenty years ago by the second
author.

\begin{con}\cite{Ow}\label{iso}
No $3$-connected minimally rigid graph is radically solvable.
\end{con}

We have verified that the smallest 3-connected non-planar minimally rigid
graph, $K_{3,3}$, is not radically solvable using a similar proof technique
to that used for the prism, or doublet, graph in
\cite[Theorem 8.4]{OP}.

\section{2-connected graphs}

Every rigid graph is 2-connected but not necessarily 3-connected. We
show in this section that the problem of deciding whether a
minimally rigid graph is radically, respectively quadratically,
solvable can be reduced to the special case of 3-connected minimally
rigid graphs. We obtain similar reduction results for arbitrary
rigid graphs but in this case the reduction to 3-connected graphs is
not complete.




\begin{thm}\label{2sep}
Let $H_1=(V_1,E_1)$ and $H_2=(V_2,E_2)$ be graphs with $V_1\cap
V_2=\{v_1,v_2\}$ and $E_1\cap E_2=\emptyset$. Let $G=H_1\cup H_2$
and suppose that $G$ is rigid.
\\
(a) Suppose that $H_1,H_2$ are both rigid. Then $G$ is radically,
respectively quadratically, solvable if and only if $H_1+v_1v_2,
H_2+v_1v_2$ are both radically, respectively
quadratically, solvable.\\
(b) Suppose that $H_1$ is not rigid. Then $H_1+v_1v_2$ and $H_2$ are
both rigid. Furthermore:\\
(i) if $H_1+v_1v_2$ and $H_2$ are both radically, respectively
quadratically,  solvable then $G$ is radically, respectively quadratically, solvable;\\
(ii) if $G$ is radically, respectively quadratically, solvable then
$H_1+v_1v_2$ and $H_2+v_1v_2$ are both
radically, respectively quadratically, solvable.\\
(iii) if $G$ is radically, respectively quadratically, solvable and
$H_1+v_1v_2$ is minimally rigid, then  $H_1+v_1v_2$ and $H_2$ are
both radically, respectively quadratically, solvable.
\end{thm}
\bproof Choose a quasi-generic realisation $(G,p)$ of $G$ with
$p(v_1)=(0,0)$ and $p(v_2)=(0,y)$ for some $y\in \complex$.

\medskip\noindent
(a) Suppose that $H_1+v_1v_2$ and  $H_2+v_1v_2$ are both radically,
respectively quadratically, solvable. Then $\rat(p|_{V_i})$ is
a radically, respectively quadratically, solvable extension of
$\rat(d_{H_i+v_1v_2}(p))$ for $i=1,2$. It follows that $\rat(p)$ is
a radically, respectively quadratically, solvable extension of
$\rat(d_{G+v_1v_2}(p))$. Since $H_1,H_2$ are both rigid, \cite[Lemma
8.2]{JO} implies that $v_1$ and $v_2$ are globally linked in
$(G,p)$. By Lemma \ref{globally_b}, $d_p(v_1,v_2)\in \rat(d_G(p))$
and hence $\rat(d_{G+v_1v_2}(p))=\rat(d_G(p))$. Thus $\rat(p)$ is
a radically, respectively quadratically, solvable extension of
$\rat(d_G(p))$ and $G$ is radically, respectively quadratically,
solvable.

Suppose on the other hand that $G$ is radically, respectively
quadratically, solvable. Then $\rat(p)$ is
a radically, respectively quadratically, solvable extension of
$\rat(d_{G}(p))$. Since $(G,p)$ is quasi-generic, we may apply Lemma
\ref{sep} with $f=d_{H_1}(p)$, $g=d_{H_2}(p)$,
$X=p|_{V_1\setminus\{v_1,v_2\}}$, $Y=y$, and
$Z=p|_{V_2\setminus\{v_1,v_2\}}$ to deduce that $\rat(p|_{V_1})$ is
a radically, respectively quadratically, solvable
extension of $\rat(d_{H_1}(p),y)$. Thus $\rat(p|_{V_1})$ is
a radically, respectively quadratically, solvable extension of
$\rat(d_{H_1+v_1v_2}(p))$ and so $H_1+v_1v_2$ is radically,
respectively quadratically, solvable. By symmetry, $H_2+v_1v_2$ is
also radically, respectively quadratically, solvable.

\medskip\noindent
(b) The fact that $H_1+v_1v_2$ and  $H_2$ are both rigid follows
from \cite[Lemma 8.5]{JO}.

Suppose that $H_1+v_1v_2$ and $H_2$ are both radically, respectively
quadratically, solvable. Then $\rat(p|_{V_2})$ is
a radically, respectively quadratically, solvable extension of
$\rat(d_{H_2}(p))$. We also have $\rat(p|_{V_1})$ is
a radically, respectively quadratically, solvable extension of
$\rat(d_{H_1+v_1v_2}(p))$. Since $y\in \rat(p|_{V_2})$, we have
$\rat(p|_{V_1},p|_{V_2}))$  is 
a radically, respectively quadratically, solvable extension of
$\rat(d_{H_1}(p),p|_{V_2})$. Thus $\rat(p)$ is
a radically, respectively quadratically, solvable extension of
$\rat(d_G(p))$ and $G$ is radically, respectively quadratically,
solvable. Hence (i) holds.

Suppose on the other hand that that $G$ is radically, respectively
quadratically, solvable. Then $\rat(p)$ is
a radically, respectively quadratically, solvable extension of
$\rat(d_{G}(p))$. We may apply the argument used in the second part
of the proof of (a) to deduce that $\rat(p|_{V_1})$ is
a radically, respectively
quadratically, solvable extension of $\rat(d_{H_1+v_1v_2}(p))$, and
$\rat(p|_{V_2})$ is
a radically, respectively quadratically, solvable extension of
$\rat(d_{H_2+v_1v_2}(p))$. Hence (ii) holds.

To prove (iii) we need to show that $y$ belongs to a radical,
respectively quadratic, extension of $\rat(d_{H_2}(p))$ when
$H_1+v_1v_2$ is minimally rigid. In this case \cite[Lemma 5.6]{JO}
implies that $X=d_{H_1}(p)$ is algebraically independent over
$\rat(d_{H_2+v_1v_2}(p))$. Let $K=\rat(d_{H_2}(p))$ and $L=K(y)$.
Since $G$ is radically, respectively quadratically, solvable,
$L(X):K(X)$ is radically, respectively quadratically, solvable.
Since $X$ is algebraically independent over $L$, Lemma \ref{indext}
implies that $L:K$ is radically, respectively quadratically,
solvable. Part (iii) now follows since $y\in L$.
\eproof

We do not know whether the hypothesis that $H_1+v_1v_2$ is minimally
rigid can be removed from Theorem \ref{2sep}(b)(iii). The difficulty
in extending the above proof when $H_1+v_1v_2$ is not minimally
rigid is that $d_{H_1}(p)$ will not be algebraically independent
over $\rat(d_{H_2+v_1v_2}(p))$. So it is conceivable that
$\rat(d_{H_1}(p))$ may contain algebraic numbers which enable $y$ to
belong to a radical extension of $\rat(d_{G}(p))$ but not to a
radical extension of $\rat(d_{H_2}(p))$. On the other hand, we will
see in the next section that we can side step this problem and still
obtain a characterization of radically solvable rigid graphs if
Conjecture 6.6 is true. We will accomplish this by only considering
certain separations $(H_1,H_2)$ of $G$ and applying the following
result.


\begin{cor}\label{triangle}
Let $H_i=(V_i,E_i)$ be graphs with $V_i\cap V_j=\{v_k\}$ and
$V_1\cap V_2\cap V_3=\emptyset=E_i\cap E_j$ for all
$\{i,j,k\}=\{1,2,3\}$. Let $G=H_1\cup H_2\cup H_3$ and suppose that
$G$ is rigid. Then $H_1,H_2,H_3$ are rigid. Furthermore, $G$ is
radically, respectively quadratically, solvable if and only if
$H_1,H_2,H_3$ are radically, respectively quadratically, solvable.
\end{cor}
\bproof Since $G=(H_1\cup H_2)\cup H_3$ is rigid and $H_1\cup H_2$
is not rigid, Theorem \ref{2sep}(b) implies that $H_3$ is rigid. We
may now use symmetry to deduce that $H_1,H_2$ are also rigid.

Suppose $G$ is radically, respectively quadratically, solvable. By
Theorem \ref{2sep}(b)(ii), $(H_1\cup H_2)+v_1v_2$ is radically,
respectively quadratically, solvable. Since $(H_1\cup
H_2)+v_1v_2=H_1\cup (H_2+v_2+v_1v_2)$ we may again use Theorem
\ref{2sep}(b)(ii) to deduce that $H_2+v_2+v_2v_3+v_1v_2$ is
radically, respectively quadratically, solvable. We can now express
$H_2+v_2+v_2v_3+v_1v_2$ as $(K_3-v_1v_3)\cup H_2$ where
$V(K_3)=\{v_1,v_2,v_3\}$. Since $K_3$ is minimally rigid, we may
apply Theorem \ref{2sep}(b)(iii) to deduce that $H_2$ is radically,
respectively quadratically, solvable. By symmetry $H_1,H_3$ are also
radically, respectively quadratically, solvable.

Suppose on the other hand that $H_1,H_2,H_3$ are radically,
respectively quadratically, solvable. Let $K_3$ be a complete graph
with $V(K_3)=\{v_1,v_2,v_3\}$. Then $K_3$ is quadratically solvable,
so by Theorem \ref{2sep}(b)(i), $F_1=(K_3-v_1v_3)\cup H_2$ is
radically, respectively quadratically, solvable. We may now apply
Theorem \ref{2sep}(b)(i) to $F_2=(F_1-v_2v_3)\cup H_1$ to deduce
that $F_2$ is radically, respectively quadratically, solvable.
Finally we apply Theorem \ref{2sep}(b)(i) to $G=(F_2-v_1v_2)\cup
H_3$ to deduce that $G$ is radically, respectively quadratically,
solvable. \eproof

\section{A family of quadratically solvable graphs}

\begin{figure}
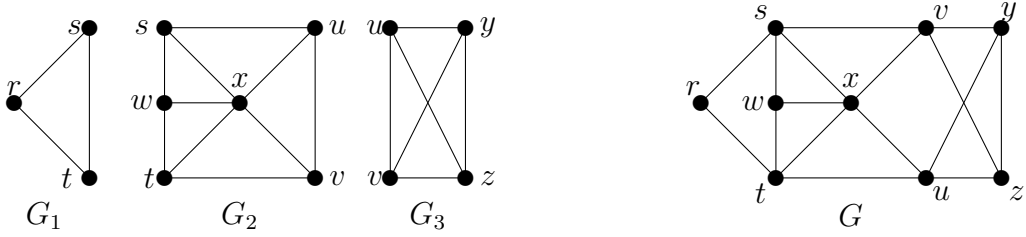

\input{radfig3.tec}
\input{radfig3a.tec}
\caption{Three globally rigid graphs $G_1,G_2,G_3$ are combined to
give a graph $G$ in $\cal F$. We first construct $G_4=(G_1-st)\cup
G_2$ using operation (b). We then construct $G=(G_4-uv)\cup
(G_3-uv)$ using operation (c).} \label{fig3}
\end{figure}

We can recursively construct a family $\cal F$ of quadratically
solvable graphs as follows. We first put all globally rigid graphs
in $\cal F$. Then, for any $G_1=(V_1,E_1)$ and $G_2=(V_2,E_2)$ in
$\cal F$ with $V_1\cap V_2=\{u,v\}$ and $|V_1|,|V_2|\geq 3$  we put:
\begin{enumerate}
\item[(a)] $G_1\cup G_2$ in $\cal F$;
\item[(b)]  $(G_1-e)\cup G_2$ in $\cal F$ if $e=uv\in E_1$;
\item[(c)]  $(G_1-e)\cup (G_2-e)$ in $\cal F$  if $e=uv\in E_1\cap
E_2$ and $G_1-e, G_2-e$ are both rigid.
\end{enumerate}
This construction is illustrated in Figure \ref{fig3}. (Note that a
recursive construction for globally rigid graphs is given in
\cite{JJ}.)

\begin{lem}\label{construct}
Every graph in $\cal F$ is rigid and quadratically solvable.
\end{lem}
\bproof Suppose $G\in \cal F$. We show that $G$ is rigid and
quadratically solvable by induction on $|E|$. If $G$ is globally
rigid then $G$ is rigid, and is quadratically solvable by Theorem
\ref{globallyrigid}. Hence we may suppose that $G$ is not globally
rigid. The definition of $\cal F$ now implies that there exist
graphs $G_1=(V_1,E_1)$ and $G_2=(V_2,E_2)$ in $\cal F$ with $V_1\cap
V_2=\{u,v\}$ and either $G= G_1\cup G_2$, or  $e=uv\in E_1$ and
$G=(G_1-e)\cup G_2$, or $e=uv\in E_1\cap E_2$, $G_1-e, G_2-e$ are
both rigid and $G=(G_1-e)\cup (G_2-e)$. By induction $G_1$ and $G_2$
are both rigid and quadratically solvable.

We first show that $G$ is rigid.  Since $G_1,G_2$ are rigid and
$|V_1\cap V_2|\geq 2$, $G_1\cup G_2$ is rigid. Furthermore, if
$e=uv\in E_1$ then $e$ is a redundant edge in $G_1\cup G_2$, so
$(G_1-e)\cup G_2$ is also rigid. Finally, if $e\in E_1\cap E_2$ and
$G_1-e$ and $G_2-e$ are both rigid then $(G_1-e)\cup (G_2-e)$ is
rigid. Hence $G$ is rigid.

It remains to show that $G$ is quadratically solvable. Since $G_1$
and $G_2$ are quadratically solvable, $G_1+uv$ and $G_2+uv$ are
quadratically solvable. Hence $G_1\cup G_2$ is quadratically
solvable by Theorem \ref{2sep}(a). Suppose that $e=uv\in E_1$ and
let $H_1=G_1-e$ and $H_2=G_2$. We can deduce that $G=H_1\cup H_2$ is
quadratically solvable by applying Theorem \ref{2sep}(a) to $H_1$
and $H_2+uv$ if $H_1$ is rigid, and by applying  Theorem
\ref{2sep}(b)(i) to $H_1$ and $H_2$ if $H_1$ is not rigid. Thus
$(G_1-e)\cup G_2$ is quadratically solvable.
  Finally we suppose that $e\in E_1\cap E_2$ and $G_1-e, G_2-e$ are both rigid.  Then
$(G_1-e)\cup (G_2-e)$ is quadratically solvable by Theorem
\ref{2sep}(a). \eproof

We can use Theorems \ref{planar} and \ref{2sep} and Lemma
\ref{construct} to characterize when a rigid planar graph is
quadratically solvable. We first need to describe a technique for
decomposing a rigid graph into `3-connected rigid pieces'. This is a
special case of a more general theory of Tutte \cite{T} which
decomposes 2-connected graphs into `3-connected pieces'.

Every 2-connected graph $G$ which is distinct from $K_3$ and is not
3-connected has a pair of edge-disjoint subgraphs
$H_1=(V_1,E_1)$ and $H_2=(V_2,E_2)$ such that $H_1\cup H_2=G$,
$|V_1\cap V_2|= 2$, and $V_1\sm V_2\neq \emptyset\neq V_2\sm V_1$.
We refer to such a pair of subgraphs $(H_1,H_2)$ as a {\em
$2$-separation} of $G$ and to the vertex set $V_1\cap V_2$ as a {\em
$2$-separator} of $G$.

Given a  rigid graph $G$ with at least three vertices, we
recursively construct the set ${\cal C}_G$ of {\em cleavage units}
of $G$ as follows. If $G$ is 3-connected or $G=K_3$ then we put
${\cal C}_G=\{G\}$. Otherwise $G$ has a 2-separation $(H_1,H_2)$,
where $V(H_1)\cap V(H_2)=\{u,v\}$. In this case $G_1=H_1+uv$ and
$G_2=H_2+uv$ are both rigid by Theorem \ref{2sep}(b), and we put
${\cal C}_G={\cal C}_{G_1}\cup{\cal C}_{G_2}$.\footnote{In order to
obtain a unique decomposition of a 2-connected graph $G$ into
cleavage units Tutte \cite{T} only considers {\em excisable
2-separations} i.e. 2-separations $(H_1,H_2)$ such that at least one
of  $H_1,H_2$ is 2-connected. When $G$ is rigid, Theorem
\ref{2sep}(b) tells us that for every 2-separation $(H_1,H_2)$, at
least one of $H_1,H_2$ will be rigid (and hence 2-connected) so all
2-separations of a rigid graph are excisable.} Note that the
cleavage units of $G$ may not be subgraphs of $G$ since $G_1$ and
$G_2$ may not be subgraphs of $G$. (We have $uv\in E(G_1)\cap
E(G_2)$ but we may not have $uv\in E(G)$. For example the cleavage
units of the graph $G$ in Figure \ref{fig3} are $G_1$, $G_2+st$ and
$G_3$, and none of these are subgraphs of $G$.)

\begin{lemma}\label{cleavage}
Let $G$ be a rigid graph on at least three vertices. Then every
cleavage unit of $G$ is either equal to $K_3$ or is 3-connected and
rigid. Furthermore, if $G$ is radically, respectively quadratically,
solvable, then every cleavage unit of $G$ is radically, respectively
quadratically, solvable.
\end{lemma}
\bproof If $G$ itself is $K_3$ or is 3-connected then the lemma is
trivially true. Hence we may suppose that $G$ has a 2-separation
$(H_1,H_2)$, where $V(H_1)\cap V(H_2)=\{u,v\}$. Theorem \ref{2sep}
implies that $H_1+uv, H_2+uv$ are both rigid, and are radically,
respectively quadratically, solvable if $G$ is radically,
respectively quadratically, solvable. The lemma now follows by
induction on $|V(G)|$ using the fact that ${\cal C}_G={\cal
C}_{H_1+uv}\cup{\cal C}_{H_2+uv}$.
\eproof

We can now obtain our promised characterization of quadratic
solvability for rigid planar graphs.

\begin{theorem}\label{planar2con}
Let $G$ be a  rigid planar graph. Then the following statements are
equivalent.\\
(a) $G$ is quadratically solvable. \\
(b) $G$ is radically solvable.\\
(c) $G$ belongs to $\cal F$.
\end{theorem}
\bproof We have (c) implies (a) by Lemma \ref{construct}, and (a)
implies (b) by definition. It remains to show that (b) implies (c).
We proceed by contradiction. Suppose there exists a radically
solvable rigid planar graph $G$ such that $G\not\in \cal F$. We may
assume that $G$ is chosen to have as few vertices as possible (and
hence every radically solvable rigid planar graph with fewer
vertices than $G$ belongs to $\cal F$). Since $G\not\in \cal F$,
$G\neq K_2,K_3$. If $G$ were 3-connected then $G$ would be globally
rigid by Theorem \ref{planar} and hence we would have $G\in \cal F$.
Thus $G$ is not 3-connected and we may choose a 2-separation
$(H_1,H_2)$ of $G$, where $V(H_1)\cap V(H_2)=\{u,v\}$. By Theorem
\ref{2sep}, $H_1+uv,H_2+uv$ are both rigid and radically solvable.
Since they are also planar and have fewer vertices than $G$ we have
$H_1+uv,H_2+uv\in \cal F$. If $uv\in E(G)$ then $G=(H_1+uv)\cup
(H_2+uv)\in \cal F$ by operation (a) in the definition of $\cal F$.
Hence $uv\not\in E(G)$. If $H_1,H_2$ are both rigid then  $G=H_1\cup
H_2\in \cal F$  by operation (c) in the definition of $\cal F$.
Thus, for every 2-separator $\{u,v\}$ of $G$, $uv\not\in E(G)$, and
for every 2-separation $(H_1,H_2)$ of $G$, one of $H_1$ and $H_2$ is
not rigid.

We now modify our choice of the 2-separation $(H_1,H_2)$ if
necessary so that $H_1$ is not rigid and, subject to this condition,
$H_1$ has as few vertices as possible.

\begin{figure}
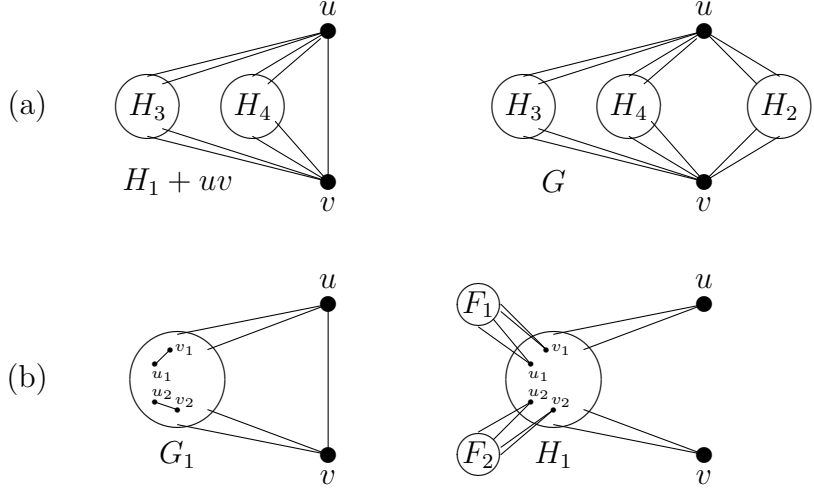

\input{radfig4a.tec}
\input{radfig4bnew.tec}
\caption{Proof of Claim \ref{c1}: (a) the case when there are two
distinct cleavage units $G_3,G_4$ of $G$ with $\{u,v\}\subset
V(G_i)\subseteq V(H_1)$ for $i=3,4$; (b) the case when $G_1\neq
K_3$.} \label{fig4}
\end{figure}

\begin{claim}\label{c1} There exists a unique cleavage unit $G_1$ of $G$ with
$\{u,v\}\subset V(G_1)\subseteq V(H_1)$. In addition we have
$G_1=K_3$.
\end{claim}
\bproof Suppose that there are two distinct cleavage units $G_3,G_4$
of $G$ with $\{u,v\}\subset V(G_i)\subseteq V(H_1)$ for $i=3,4$.
Then $H_1+uv$ has a 2-separation $(H_3,H_4)$ with $uv \in E(H_4)$,
$V(H_3)\cap V(H_4)=\{u,v\}$ and $V(G_i)\subseteq V(H_i)$ for
$i=3,4$, see Figure \ref{fig4}(a). Since $H_1=H_3\cup (H_4-uv)$ is
not rigid, $H_3$ is not rigid. Thus $(H_3, H_2\cup (H_4-uv))$ is a
2-separation of $G$ in which $H_3$ is not rigid and has fewer
vertices than $H_1$. This contradicts the choice of $(H_1,H_2)$.
Hence there is a unique cleavage unit $G_1$ of $G$ with
$\{u,v\}\subset V(G_1)\subseteq V(H_1)$.

Suppose $G_1\neq K_3$. Then $G_1$ is 3-connected and radically
solvable by Lemma \ref{cleavage}.
Since $G_1$ is planar, Theorem \ref{planar} now implies that $G_1$
is redundantly rigid and hence that $G_1-uv$ is rigid. Let
$\{u_i,v_i\}$, $1\leq i \leq m$, be the 2-separators of $H_1+uv$
with $\{u_i,v_i\}\subset V(G_1)$. Then $u_iv_i\in E(G_1)$ for
$1\leq i\leq m$, see Figure \ref{fig4}(b). For each
$1\leq i \leq m$ we may choose a  2-separation $(F_i,F_i')$ of
$H_1+uv$ with $V(G_1)\subset V(F_i')$. Then $(F_i,(F_i'-uv)\cup
H_2)$ is a 2-separation of $G$. The choice of $H_1$ and the fact
that $F_i$ is properly contained in $H_1$ now implies that $F_i$ is
rigid for all $1\leq i\leq m$. Since $G_1-uv$ is rigid, this implies that
$$H_1=[(G_1-uv)-\{u_iv_i\,:\,1\leq i\leq m\}]\cup
\bigcup_{i=1}^mF_i$$ is rigid. This contradicts the choice of $H_1$.
Thus $G_1=K_3$. \eproof

We can now complete the proof of the theorem. Since $G_1=K_3$ we can
express $G$ as $G=H_1'\cup H_1''\cup H_2$ where $H_1'\cup
H_1''=H_1$, $V(H_1')\cap V(H_2)=\{u\}$, $V(H_1'')\cap V(H_2)=\{v\}$,
$V(H_1')\cap V(H_1'')=\{w\}$ for some $w\in V(H_1)\sm \{u,v\}$, and
$H_1',H_1'',H_2$ are pairwise edge-disjoint. Corollary
\ref{triangle} now implies that $H_1',H_1'',H_2$ are rigid and
radically solvable. Since they are planar and have fewer vertices
than $G$, we have $H_1',H_1'',H_2\in \cal F$. Since $G$ can be
obtained from $K_3,H_1',H_1'',H_2$ by applying operation (b) in the
definition of $\cal F$ at most three times, we have $G\in \cal F$.
This contradicts the choice of $G$. \eproof

Since the operations (a), (b) and (c) used in the construction of
$\cal F$ preserve planarity, Theorem \ref{planar2con} implies that
the family of quadratically solvable planar graphs can be
constructed recursively from the family of globally rigid planar
graphs by applying operations (a), (b) and (c).

We conjecture that the planarity hypothesis can be removed from
Theorem  \ref{planar2con}.

\begin{con}\label{planechar}
Let $G$ be a  rigid graph. Then the following statements are
equivalent.\\
(a) $G$ is quadratically solvable. \\
(b) $G$ is radically solvable.\\
(c) $G$ belongs to $\cal F$.
\end{con}

The proof technique of Theorem \ref{planar2con} can be used to show
that Conjecture \ref{planechar} is equivalent to Conjecture
\ref{gen}, and hence is also equivalent to Conjecture \ref{iso}.

\medskip

The constructions and some of the results of this section extend
earlier work of the second author for minimally rigid graphs which
is implicitly given in \cite{Ow}, and explicitly stated in
\cite[Theorem 3.2]{OP}. He recursively constructs a subfamily,
$\Fiso$, of $\cal F$ as follows. He first puts $K_3$ in $\Fiso$.
Then, for any $G_1=(V_1,E_1)$ and $G_2=(V_2,E_2)$ in $\Fiso$ with
$V_1\cap V_2=\{u,v\}$ and $e=uv\in E_1$ he puts $(G_1-e)\cup G_2$ in
$\Fiso$. He shows that every graph in $\Fiso$ is minimally rigid and
quadratically solvable and conjectures that every radically solvable
minimally rigid graph belongs to $\Fiso$.\footnote{Since the radical
solvability of a graph is preserved by the addition of edges, it is
tempting to also conjecture that a graph is radically solvable if
and only if it has a spanning subgraph in $\Fiso$. This is not the
case however. The complete bipartite graph $K_{3,4}$ is globally
rigid and hence quadratically solvable,  but for each edge $e$,
$K_{3,4}-e$ is minimally rigid and does not belong to $\Fiso$.
In addition, we can use Theorem \ref{2sep}(b)(iii) and the fact that
$K_{3,3}$ is not
radically solvable to deduce that $K_{3,4}-e$ is not radically solvable, so
$K_{3,4}$ is `minimally radically solvable' but not minimally rigid.
}

\medskip
\noindent {\bf Acknowledgement} We would like to thank John Bray for
helpful comments. We would also like to thank the Fields Institute
for support during its 2011 thematic programme on Discrete Geometry
and Applications.

\end{document}